\documentclass[11pt]{article}

\usepackage{amssymb,latexsym}
\usepackage{amsmath,amscd}
\usepackage{theorem}
\usepackage{url}
\usepackage{bm} 
\usepackage{graphics}

\title{\bf Index-stable compact $\bm{p}$-adic analytic  groups}

\author{
        Francesco Noseda
        \\[0.1cm]
        Ilir Snopce\thanks{Supported by CNPq
        and FAPERJ.}
        \\[0.2cm]
        {\footnotesize \textit{with an appendix by}}  Jean-Pierre Serre\\      
}

\date{}

\newcommand{\bb}[1]{\mathbb{#1}}
\newcommand{\cl}[1]{\mathcal{#1}}

\newcommand{\mr}[1]{\mathrm{#1}}

\newcommand{\les}{\leqslant}
\newcommand{\ges}{\geqslant}

\newcommand{\ep}{\hfill $\square$} 

\def\Q{\mathbb Q}

\newtheorem{lemma} {Lemma} 

\newtheorem{theorem} [lemma] {Theorem}
\newtheorem{corollary} [lemma] {Corollary}

\theorembodyfont{\rm} 

\newtheorem{remark}[lemma]{Remark}

\newtheorem{question}[lemma]{Question}

\numberwithin{equation}{section}

\theorembodyfont{\it}

\usepackage{array}
\usepackage[colorlinks]{hyperref}

\begin{document} 

\maketitle

\begin{abstract}
A profinite group is index-stable
if any two isomorphic open subgroups have the same index.
Let $p$ be a prime, and let $G$ be a compact $p$-adic analytic group with 
associated $\mathbb{Q}_p$-Lie algebra $\cl{L}(G)$. 
We prove that  $G$ is index-stable
whenever $\cl{L}(G)$ is semisimple.
In particular, a just-infinite  compact $p$-adic analytic group is index-stable if and only if it is not virtually abelian. 
Within the category of compact $p$-adic analytic groups, this gives a positive answer to a question of C. Reid.

In the Appendix, J-P. Serre proves that $G$ is index-stable if and only if
the determinant of any automorphism of $\cl{L}(G)$ has $p$-adic norm 1. 
\end{abstract}
{
\let\thefootnote\relax
\footnotetext{\textit{Mathematics Subject Classification (2020): }
Primary 20E18, 22E20; Secondary 22E60.}
\footnotetext{\textit{Key words:} index-stable group,
just-infinite profinite group, 
$p$-adic analytic group, pro-$p$ group, $p$-adic Lie lattice, 
commensurator.}
}


\section*{Introduction}

Throughout, let $p$ be a prime. A pro-$p$ groug $G$ is said to be 
\emph{powerful} if $p \ges 3$ and
$[G,G]\les G^p$, or $p=2$ and $[G,G]\les G^4$.   
A finitely generated torsion-free powerful pro-$p$ group is called \emph{uniform}.   
In his seminal paper \emph{Groupes analytiques $p$-adiques} \cite{Laz65}, M. Lazard obtained the following algebraic characterization of $p$-adic analytic groups:
a topological group is $p$-adic analytic if and only if it contains an open uniform pro-$p$ subgroup (see \cite[Theorems 8.1 and 8.18]{DixAnaProP}). 
With every uniform pro-$p$ group $U$ one can naturally associate a $\mathbb{Z}_p$-Lie lattice $L_U$. 
The $ \mathbb{Q}_p$-Lie algebra associated with a compact $p$-adic analytic group $G$ is defined as $\mathcal{L}(G): = L_U \otimes_{\bb{Z}_p} \mathbb{Q}_p$, where $U$ is an open uniform pro-$p$ subgroup of $G$.

\medskip

The following theorem is the main result of this paper.

\begin{theorem}\label{index-group}
Let $G$ be a  compact $p$-adic analytic group, let $\mathcal{L}(G)$ be the
$\mathbb{Q}_p$-Lie algebra associated with $G$, and assume that $\mathcal{L}(G)$
is semisimple.
If  $H$ and $K$ are two isomorphic closed subgroups  of $G$ then $|G:H|= |G:K|$.
\end{theorem}

\noindent
\textbf{Remark}\hspace{3pt}
Soon after a preliminary version of this paper was published on the arXiv, the 
authors received a letter from Jean-Pierre Serre with a different proof
of a more general version of Theorem \ref{index-group};
see Theorem \ref{thserre}.
We thank Professor Serre for kindly agreeing to include his letter as an appendix to this paper.
\\

A profinite group is said to be \emph{index-unstable} 
if it contains a pair of isomorphic open subgroups of different indices;
otherwise, it is said to be \emph{index-stable}. 
This definition was introduced by C. Reid in \cite{ReiJust-Infinite}, 
where he also raised the following question, which is still open.

\begin{question}\label{question1}
Let $G$ be a (hereditarily) just-infinite profinite group which is index-unstable. 
Is $G$ necessarily virtually abelian?
\end{question}

\smallskip

Recall that a profinite group $G$ is said to be \emph{just-infinite} if it is infinite, and every non-trivial closed normal subgroup of $G$ is of finite index. 
A just-infinite profinite group $G$ is  \emph{hereditarily just-infinite} if every open subgroup of $G$ is just-infinite.

\smallskip

The following corollary gives a positive answer to Question \ref{question1} within the category of $p$-adic analytic groups.

\begin{corollary}\label{index-stable-Reid}
Let $G$ be a  just-infinite compact $p$-adic analytic group. 
Then $G$ is index-stable if and only if it is not virtually abelian. 
\end{corollary}

\smallskip

Question \ref{question1} was motivated by the study of commensurators of profinite groups, 
in the spirit of Y. Barnea, M. Ershov and T. Weigel \cite{BEWAbs-Com}.  Let us recall briefly some relevant definitions.
A virtual automorphism of a profinite group $G$ is  an isomorphism between two open subgroups of $G$. Two virtual automorphisms of $G$ are said to be equivalent if they agree on some open subgroup of $G$. 
Given two virtual automorphisms, one can define their composition up to equivalence by composing suitable equivalence-class representatives.
 With respect to this composition, the  equivalence classes of virtual automorphisms of $G$ form a group $\textrm{Comm}(G)$ which is called the \emph{commensurator} of $G$.
 In \cite{BEWAbs-Com}, a topology is defined on $\textrm{Comm}(G)$, called the strong topology, so that it becomes a 
  topological group $\textrm{Comm}(G)_S$.
The \emph{virtual center} of $G$ is defined to be 
$V\!Z(G) = \{ x \in G :\, \textrm{Cent}_G(x) \textrm{ is open in }  G \}$. 
If $V\!Z(G) =\{1\}$ then $\textrm{Comm}(G)_S$ is locally compact;
in this case, $\textrm{Comm}(G)_S$ is said to be \emph{unimodular}
if its modular function is trivial, or equivalently, if the Haar measure is both left and right invariant.
 
Suppose that $U$ and $V$ are isomorphic open subgroups of $G$. 
 Given an isomorphism $\theta: U \to V$, let $\iota(\theta) = \frac{\lvert G:U \rvert}{\lvert G:V \rvert}$; note that this is invariant under equivalence.
  Hence, given $\phi \in \textrm{Comm}(G)$, one can define $\iota(\phi)$ as $\iota(\theta)$ for any  $\theta$ that represents $\phi$. 
  This yields a function $\iota : \textrm{Comm}(G) \to \mathbb{Q}_{+}$, where $\mathbb{Q}_{+}$ is the multiplicative group of positive rationals; this function is called the \emph{index ratio}.
If the virtual center of $G$ is trivial then the index ratio turns out to be the
modular function for the strong topology. 
Moreover, it is not difficult to see that if $G$ is a just-infinite profinite group then $V\!Z(G) =\{1\}$ if and only if $G$ is not virtually abelian.
  
\bigskip

The following result is a direct consequence of Corollary \ref{index-stable-Reid}  
and the above discussion.

 \begin{corollary}
 Let $G$ be a just-infinite compact $p$-adic analytic group.
  If $G$ is not virtually abelian then $\mathrm{Comm}(G)_S$ is unimodular.
 \end{corollary}

Solvable just-infinite $p$-adic analytic pro-$p$ groups are irreducible 
$p$-adic space groups; in particular, they are virtually abelian (cf.  \cite{KLP-Linear}). 
Hence, they contain many pairs of open subgroups of different indices that are isomorphic to each other. In contrast, 
by \cite[Lemma III.11]{KLP-Linear},  an \emph{unsolvable} just-infinite $p$-adic analytic pro-$p$ group $G$
has the remarkable property of not being isomorphic to any of its proper closed subgroups.
C. Reid generalized this result to all profinite groups by proving that a just-infinte profinite group $G$  that contains an open proper subgroup  isomorphic to $G$ is virtually abelian (see \cite[Theorem~E]{ReiJust-Infinite}). 
The following corollary shows that within the category of compact $p$-adic analytic groups a much stronger result holds.

\begin{corollary}\label{just-infinite-pro-p}
Let $G$ be an  unsolvable just-infinite $p$-adic analytic pro-$p$ group. 
Then $G$ is index-stable.
\end{corollary}

Open pro-$p$ subgroups  of  $p$-adic Chevalley groups form a rich source  of unsolvable (hereditarily) just-infinite $p$-adic analytic pro-$p$ groups; for $n\ges 2$, 
the first congruence subgroup $\textrm{SL}_n^1(\mathbb{Z}_p)$ of  $\mathrm{SL}_n(\mathbb{Z}_p)$ and its open subgroups are typical examples of such groups.
More generally, given a simple finite dimensional $\mathbb{Q}_p$-Lie algebra $\mathcal{L}$, any open pro-$p$
 subgroup of  $\textrm{Aut}(\mathcal{L})$ is an unsolvable (hereditarily) just-infinite $p$-adic analytic pro-$p$ group (see  \cite[Proposition III.9]{KLP-Linear}). Moreover, 
 given  an unsolvable just-infinite $p$-adic analytic pro-$p$ group $G$, there is a semisimple $\mathbb{Q}_p$-Lie algebra $\mathcal{L}$ such that
 $G$ is an open subgroup of  $\textrm{Aut}(\mathcal{L})$ (cf.  \cite[Section III.9]{KLP-Linear}).

\medskip

By  \cite[Proposition 6.1]{GSKpsdimJGT}, every solvable 
just-infinite pro-$p$ group other than $ \mathbb{Z}_p$ has torsion. 
Thus, a non-procyclic  torsion-free just-infinite $p$-adic analytic pro-$p$ group
must be unsolvable. Hence, we can deduce the following result, which 
is the correct formulation of \cite[Conjecture~2.10]{NSSelf}.

\begin{corollary}\label{just-torsion-free}
Let $G$ be a non-procyclic torsion-free just-infinite $p$-adic analytic pro-$p$ group. 
Then $G$ is index-stable.
\end{corollary}

Next, 
we observe that Theorem \ref{index-group} relies on its Lie-algebra counterpart,
which is an interesting result on its own.
By definition, a $\bb{Z}_p$-Lie lattice is a $\bb{Z}_p$-Lie algebra the
underlying module of which is finitely generated and free.
A $\bb{Z}_p$-Lie lattice $L$ is said to be \emph{index-stable} if for any pair 
$M$ and $N$ of
isomorphic finite-index subalgebras of $L$ we have $|L:M|=|L:N|$.
 
 \begin{theorem}\label{index-Lie}
Let $L$ be a $\bb{Z}_p$-Lie lattice, and assume that $L\otimes_{\bb{Z}_p} \bb{Q}_p$
is semisimple as $\bb{Q}_p$-Lie algebra. Then $L$ is index-stable.
\end{theorem}

An $n$-dimensional $\mathbb{Z}_p$-Lie lattice  $L$ is said to be \emph{just infinite} if 
every non-zero ideal of $L$ has dimension $n$. The following corollary proves
(the correct formulation of) \cite[Conjecture~2.9]{NSSelf}.

\begin{corollary}\label{index-stable-Lattice}
Let $L$ be a just-infinite $\bb{Z}_p$-Lie lattice, and assume
that $\mathrm{dim}\,L >1$. 
Then $L$ is index-stable.
\end{corollary}

\noindent
We remark that Corollary \ref{just-infinite-pro-p} and Corollary \ref{index-stable-Lattice} 
have applications to the study of self-similar actions
of hereditarily just-infinite $p$-adic analytic pro-$p$ groups using the language of virtual endomorphisms 
(cf.  \cite{NSSelf}).

\medskip

We close the introduction by stating the theorem proved by Serre in the Appendix.

\begin{theorem} {\normalfont\textbf{(Serre)}}\label{thserre}
Let $G$ be a compact $p$-adic analytic group, and let $\cl{L}(G)$ be the $\bb{Q}_p$-Lie algebra
associated with $G$. Then $G$ is index-stable if and only if
all automorphisms $s$ of $\cl{L}(G)$ satisfy $|\mr{det}(s)| = 1$, where
$|\cdot |$ is the $p$-adic norm of $\bb{Q}_p$.
\end{theorem}

\begin{remark}
The proof of Theorem \ref{thserre} (see the Appendix) relies on $p$-adic integration,
for which the reader may consult, for instance,
\cite[Section 7.4]{Igu00}.
Now,  denote by (*) the condition on $\cl{L}(G)$ in the statement of Theorem \ref{thserre}. 
We note that there exist $\bb{Q}_p$-Lie algebras that satisfy condition (*) but are not semisimple.
J. L. Dyer
constructed
a nilpotent Lie algebra of dimension 9 and nilpotency class 6
with the property that its group of automorphisms is unipotent \cite{Dyer70}.
In general, by a result of G. Leger and E. Luks \cite[Theorem (*)]{Leger72},
if the automorphism group $\mr{Aut}(\cl{L})$ of a nilpotent Lie algebra $\cl{L}$ of dim $>1$
is nilpotent then it is unipotent.
J. Dixmier and W. G. Lister constructed an
example of a nilpotent Lie algebra $\cl{M}$ of dimension 8 
and nilpotency class 3 such that $\mr{Aut}(\cl{M})$ is not nilpotent but the derivation algebra 
$\mr{Der}(\cl{M})$ is nilpotent \cite{Dixmier57}.
The latter condition implies that any automorphism of $\cl{M}$
has eigenvalues that are roots of unity.
Clearly, all of these Lie algebras satisfy condition (*).
\end{remark}

\begin{corollary}
Let $G$ be a compact $p$-adic analytic group of dimension $>1$, let $\cl{L}$ be the $\bb{Q}_p$-Lie algebra
associated with $G$, and assume that $\cl{L}$ is nilpotent. 
If $\mr{Aut}(\cl{L})$ is nilpotent  or $\mr{Der}(\cl{L})$
is nilpotent then $G$ is index-stable. 
\end{corollary}

  \section{Proofs of the main results}
  
  In this section we prove  Theorem \ref{index-group}, Theorem  \ref{index-Lie}, Corollary \ref{index-stable-Reid} and Corollary \ref{index-stable-Lattice}. 
  
   \vspace{5mm}

\noindent
\textbf{Proof of Theorem  \ref{index-Lie}.}
Let $\mathcal{L}: = L \otimes_{\bb{Z}_p} \mathbb{Q}_p$   
and fix a $\mathbb{Z}_p$-basis $\{a_1, ..., a_d \}$ of $L$. 
Then, clearly, $\{a_1, ..., a_d \}$ is a $\mathbb{Q}_p$-basis of $\mathcal{L}.$
Let $\kappa: \mathcal{L} \times \mathcal{L} \to \mathbb{Q}_p$, defined by $\kappa(x, y) = \textrm{tr}(\textrm{ad}(x)\textrm{ad}(y))$, 
  be the Killing form of $\mathcal{L}$. 
  Denote by $A$ the matrix representing 
  $\kappa$ with respect to the given basis; in other words, 
  $A = (A_{ij})$, where $A_{ij} = \kappa (a_i, a_j)$ for $1 \les i, j \les d.$  Since $\mathcal{L}$ is a semisimple Lie algebra over a field of characteristic 0,
   by Cartan's criterion \cite[Page 69]{JacLieA}
    we have that $\kappa$ is a non-degenerate bilinear form; in particular, $\det(A) \neq 0$.
Let  $\{m_1, ..., m_d \}$ be a $\mathbb{Z}_p$-basis of $M$, and let $\varphi: M \to N$ be an isomorphism of $\bb{Z}_p$-Lie lattices. 
Then $\{\varphi(m_1), ..., \varphi(m_d) \}$ is a $\mathbb{Z}_p$-basis of $N$. 
Let $B= (B_{ij})$ and $C= (C_{ij})$ be the change-of-basis matrices
defined by  $m_j = \sum_i B_{ij} a_i$ and 
$\varphi(m_j) = \sum_i C_{ij} a_i$. 
Note that $[L:M]= p^{v_p(\textrm{det}(B))}$ and 
$[L:N] = p^{v_p(\textrm{det}(C))}$,
where $v_p$ is the $p$-adic valuation
(one way to prove this claim is to recall that
there exist a basis $\{b_i\}$ of $L$ and non-negative integers
$\{k_i\}$ such that $\{p^{k_i}b_i\}$ is a basis of $M$;
similarly for $N$; cf. \cite[Lemma 10.7.2]{Grillet99}).
The matrices of $\kappa$ with respect to the bases $\{m_i\}$ and $\{\varphi(m_i)\}$
are $B^TAB$ and $C^TAC$, respectively.
Since the automorphism of $\cl{L}$ induced by $\varphi$ 
preserves the Killing form, $B^TAB= C^TAC$.
Taking the determinant on both sides, and recalling that $\mr{det}(A)\neq 0$,
we see that $v_p(\mr{det}(B))=v_p(\mr{det}(C))$, and the theorem follows. 
\ep

\vspace{5mm}

\noindent
\textbf{Proof of Corollary \ref{index-stable-Lattice}.}
Since $L$ is a just-infinite $\mathbb{Z}_p$-Lie lattice of dimension greater than 1, 
it is not difficult to see that $L\otimes_{\bb{Z}_p} \bb{Q}_p$ 
is a simple $\mathbb{Q}_p$-Lie algebra. Hence, the corollary follows from Theorem \ref{index-Lie}.
\ep

\vspace{5mm}

A $\mathbb{Z}_p$-Lie lattice $L$ is called \emph{powerful} if
$[L,L] \subseteq pL$ for $p$ odd, or $[L,L] \subseteq 4L$ for $p=2$. 

With a uniform pro-$p$ group $G$ one may associate a powerful $\mathbb{Z}_p$-Lie lattice $L_G$ in the following way: $G$ and $L_G$ are identified as sets, 
and the Lie operations are defined by
\begin{displaymath}
g+h=\lim_{n \to \infty}(g^{p^n}h^{p^n})^{p^{-n}},  ~ ~ ~  [g,h]_{\mr{Lie}}=\lim_{n \to \infty}[g^{p^n},h^{p^n}]^{p^{-2n}}=\lim_{n\to \infty} (g^{-p^n}h^{-p^n}g^{p^n}h^{p^n})^{p^{-2n}} .
\end{displaymath}
On the other hand, if $L$ is a powerful $\mathbb{Z}_p$-Lie lattice, then the Campbell-Hausdorff formula induces a group structure on $L$; the resulting group is a uniform pro-$p$ group. If this construction is applied to the $\mathbb{Z}_p$-Lie Lattice $L_G$ associated with a uniform group $G$, one recovers the original group. 
Indeed, the assignment $G\mapsto L_G$ gives an isomorphism between the category of uniform pro-$p$ groups and the category of powerful $\mathbb{Z}_p$-Lie lattices (see \cite[Theorems 4.30 and 9.10]{DixAnaProP}). 
Recall that every compact $p$-adic analytic group contains an open subgroup that is a uniform pro-$p$ group  (\cite[Corollary 8.34]{DixAnaProP}). 
As we already mentioned in the introduction, the $ \mathbb{Q}_p$-Lie algebra associated 
with a compact $p$-adic analytic group $G$ is defined as 
$\mathcal{L}(G): = L_U \otimes_{\bb{Z}_p} \mathbb{Q}_p$, where $U$ is an open uniform pro-$p$ subgroup of $G$ (see \cite[Section 9.5]{DixAnaProP}). 
A key invariant of a $p$-adic analytic group $G$ is its dimension as a $p$-adic manifold,  denoted by $\textrm{dim}(G)$.  
Algebraically, $\textrm{dim}(G)$ can be described  as $d(U)$, 
where $U$ is any uniform open pro-$p$ subgroup of $G$ and $d(U)$ denotes the minimal cardinality of a topological generating set for $U$.

\vspace{5mm}

\noindent
\textbf{Proof of Theorem \ref{index-group}.}
Let $H$ and $K$ be two isomorphic closed subgroups of $G$. If one of these subgroups, say $H$,  has infinite index in $G$,  
then $\textrm{dim}(H) < \textrm{dim}(G)$. Hence $\textrm{dim}(K) < \textrm{dim}(G)$, and therefore $| G:K|$ is infinite as well. 
Since $G$ is virtually pro-$p$,  it is  not difficult to see that   $| G:K|$ and $| G:H|$  coincide as supernatural numbers.
 
Now suppose that $H$ and $K$ are open subgroups in $G$. Let $\varphi: H \to K$  be an isomorphism, and let $U$ be an open uniform subgroup of $G$; note that, since $G$ is compact, $U$ is of finite index in $G$. Let $U_H$ be a uniform open subgroup of  $H$. Then $U_K:  = \varphi(U_H)$ is a uniform open subgroup of $K$. Moreover, 
$| H:U_H | = | K:U_K |.$ 
Choose a positive integer $p^m$ such that $V_H : = {(U_H)}^{p^m}$ and  $V_K : = {(U_K)}^{p^m}$ are contained in $U$; 
by \cite[Theorem~3.6]{DixAnaProP},  $V_H$ and $V_K$ are uniform. Let $L_U$ be the $\bb{Z}_p$-Lie lattice associated with $U$. 
Then  the $\mathbb{Q}_p$-Lie algebra associated with $G$ is given by $\mathcal{L}(G) = L_U\otimes_{\bb{Z}_p} \mathbb{Q}$,
which is semisimple by assumption. 
Since $U_H$ and $U_K$ are isomorphic, we have that $V_H$ and $V_K$ are isomorphic open uniform subgroups of $U$. 
This implies that $L_{V_H}$ and $L_{V_K}$ are isomorphic $\bb{Z}_p$-Lie lattices of $L_U$ of finite index. By Theorem \ref{index-Lie},
$| L_U:L_{V_H} | = | L_U:L_{V_K} |$. Hence, by \cite[Proposition 4.31]{DixAnaProP}, $| U: V_H | = | U:V_K |$. 
Clearly, this implies that $ | G: V_H | = | G:V_K |$. Now, since  $| H:U_H | = | K:U_K |$
and  $ | U_H:V_H | = | U_K:V_K |$, it follows that  $| G:H | = | G:K |$, as desired.   
\ep

\vspace{5mm}

\noindent
\textbf{Proof of Corollary \ref{index-stable-Reid}.}
Clearly, if $G$ is virtually abelian then $G$ is index-unstable.
Suppose that $G$ is not virtually abelian.
By \cite[Theorem~3.6]{Gri-branch}, a just-infinite profinite group is either a branch profinite group or it contains an open normal subgroup 
which is isomorphic to the direct product of a finite number of copies  of some hereditarily just-infinite profinite group. Note that $G$ cannot be a branch profinite group since,  by  \cite[Theorem~8.33]{DixAnaProP}, 
$G$ is of finite rank (that is, there exists a positive integer $d$ such that every closed subgroup of $G$ can be  generated topologically by at most $d$ elements). 
Hence, $G$ contains an open normal subgroup $H$ 
which is isomorphic to the direct product of a finite number, say $k$, of copies  of some hereditarily just-infinite uniform pro-$p$ group $U$; in particular, $H$ is uniform.  Note that if $U$ is solvable  
then $U$ is isomorphic to $\mathbb{Z}_p$, and consequently $G$ is virtually abelian, 
which is a contradiction. Hence, $U$ is an unsolvable hereditarily just-infinite uniform pro-$p$ group. 
By \cite[Proposition F]{GSKpsdimJGT}, the associated $\mathbb{Q}_p$-Lie algebra 
$\mathcal{L}(U) = L_U\otimes_{\bb{Z}_p} \mathbb{Q}_p$ is simple.  Thus, the $\mathbb{Q}_p$-Lie algebra $\mathcal{L}(G) = L_H\otimes_{\bb{Z}_p} \mathbb{Q}_p$ associated with $G$ is semisimple,   since it is isomorphic to the direct sum of $k$ copies of  $\mathcal{L}(U)$.
Now Theorem  \ref{index-group} yields that $G$ is index-stable.
\ep 


\section*{Appendix - A letter from Jean-Pierre Serre}

\hspace{\parindent}Dear MM. Noseda and Snopce,\\

I have just seen your arXiv paper on ``Index-stable compact $p$-adic analytic
groups''. Your proof uses Lazard's very nice results, but in fact  these results are not 
necessary: $p$-adic integration is enough. Let me explain:
\\

Let $L$ be a finite-dimensional Lie algebra over $\bb{Q}_p$. Denote by $|x|$ the 
$p$-adic norm of $\bb{Q}_p$. Consider the following property of $L$:
\\

  (*) For every automorphism  $s$  of $L$, we have $|\mr{det}(s)| = 1$.\\

\noindent
This is true for instance if $L$ is semisimple since $\mr{det}(s)=\pm 1$, which is the case
you consider.
\\

    Assume property (*). Let $n = \dim L.$ Let $u$ be a nonzero element of $\wedge^nL^*$, where $L^*$ is the $\Q_p$-dual of $L$. Let $G$ be a compact $p$-adic analytic group with Lie algebra $L$. Then $u$ defines a right-invariant differential form $\omega_u$ on $G$ of degree $n$. The corresponding measure $\mu_u = |\omega_u|$ is a non-zero right-invariant positive measure on $G$, hence is a Haar 
    measure since $G$ is compact. Property (*) implies that $\mu_u$ is 
    invariant by every automorphism of $L$. Hence every isomorphism of $G$ onto another group $G'$ carries $\mu_u$ (for $G$) into $\mu_u$ (for $G'$). This implies that, if $G, G'$ are open subgroups of some compact $p$-adic group $G''$, then they have the same index -  as
  wanted. 
\\  

  Conversely, if (*) is not true for some  $s$, and if $G''$ is compact with Lie algebra $L$,
$s$ defines a local automorphism of  $G''$, and if  $G$  is a small enough open subgroup
of $G''$, it is transformed by  $s$  into another open subgroup  $G'$, and the ratio
$(G'':G)/(G'':G')$ is equal to $|\mr{det} (s)|$, which is  $\neq 1$.
\\
  
  Best wishes,\\

  J-P. Serre


\begin{footnotesize}

\providecommand{\bysame}{\leavevmode\hbox to3em{\hrulefill}\thinspace}
\providecommand{\MR}{\relax\ifhmode\unskip\space\fi MR }
\providecommand{\MRhref}[2]{%
  \href{http://www.ams.org/mathscinet-getitem?mr=#1}{#2}
}
\providecommand{\href}[2]{#2}

\vspace{5mm}
\newpage

\noindent
Francesco Noseda\\
Mathematics Institute\\
Federal University of Rio de Janeiro\\
Av. Athos da Silveira Ramos, 149\\
21941-909, Rio de Janeiro - RJ,
Brazil  \\
{\tt noseda@im.ufrj.br}\\

\noindent
Ilir Snopce \\
Mathematics Institute\\
Federal University of Rio de Janeiro\\
Av. Athos da Silveira Ramos, 149\\
21941-909, Rio de Janeiro - RJ,
Brazil  \\
{\tt ilir@im.ufrj.br}\\

\noindent
Jean-Pierre Serre\\ 
Collège de France\\
3, rue d'Ulm\\
75231, Paris,
France\\
{\tt serre@noos.fr}

\end{footnotesize}

\end{document}